\documentclass{article}

\usepackage{amssymb}

\newcommand{\bref}[1]{(\ref{#1})}	
\newcommand{\diagspace}{\mbox{\hspace{2em}}}	
\newcommand{\twid}[1]{\widetilde{#1}}

\newcommand{\hz}{z} 

\newcommand{\epsln}{\varepsilon}	
\newcommand{\demph}[1]{\textbf{#1}}	
\newcommand{\done}{\hfill\ensuremath{\Box}}
\newcommand{\implies}{\ \Rightarrow\ }
\newcommand{\cat}[1]{\mathcal{#1}}	
\newcommand{\fcat}[1]{\mathbf{#1}}	
\newcommand{\iso}{\,\cong\,}
\newcommand{\arrow}{\longrightarrow}
\newcommand{\isoarrow}{\stackrel{\sim}{\arrow}}

\newcommand{\csg}{*} 

\newcommand{\High}{H} 
\newcommand{\integers}{\mathbb{Z}}
\newcommand{\nat}{\mathbb{N}}   
\newcommand{\rat}{\mathbb{Q}}   
   
\newenvironment{proof}{\begin{trivlist}\item\textbf{Proof}\ }{\end{trivlist}}

\newcounter{bean}

\newtheorem{thm}{Theorem}[section]
\newtheorem{propn}[thm]{Proposition}
\newtheorem{lemma}[thm]{Lemma}
\newtheorem{cor}[thm]{Corollary}
\newtheorem{preexamples}[thm]{Examples}
\newenvironment{examples}{\begin{preexamples}\upshape}{\end{preexamples}}
\newtheorem{defn}[thm]{Definition}

\title{Objects of Categories as Complex Numbers}
\author{Marcelo Fiore\thanks{Research supported by an EPSRC Advanced
    Research Fellowship}\\
        Computer Laboratory, University of Cambridge\\
        United Kingdom\\
	Marcelo.Fiore@cl.cam.ac.uk
        \and 
        Tom Leinster\thanks{William Hodge Fellow, IH\'ES}\\
        Institut des Hautes \'Etudes Scientifiques \\
        France\\
	leinster@ihes.fr}
\date{}

\begin{document}

\maketitle

\begin{abstract}
In many everyday categories (sets, spaces, modules, \ldots) objects can be
both added and multiplied.  The arithmetic of such objects is a challenge
because there is usually no subtraction.  We prove a family of cases of
the following principle: if an arithmetic statement about the objects can
be proved by pretending that they are complex numbers, then there also
exists an honest proof.
\end{abstract}

\paragraph*{}

Consider the following absurd argument concerning planar, binary, rooted,
unlabelled trees (Blass~\cite{Blass}).  Every such tree is either the
trivial tree or consists of a pair of trees joined together at the root, so
the set $T$ of trees is isomorphic to $1 + T^2$.  Pretend that $T$ is a
complex number and solve the quadratic $T = 1 + T^2$ to find that $T$ is a
primitive sixth root of unity and so $T^6 = 1$.  Deduce that $T^6$ is a
one-element set; realize immediately that this is wrong.  Notice that $T^7
\iso T$ is, however, not obviously wrong, and conclude that it is therefore
right.  In other words, conclude that there is a bijection $T^7 \iso T$
built up out of copies of the original bijection $T \iso 1 + T^2$: a tree
is the same as seven trees.

The point of this paper is to show that `nonsense proofs' of this kind are,
actually, valid.  Our main result is approximately this:
\begin{quote} \itshape
  Let $p$, $q_1$ and $q_2$ be polynomials over $\nat$.  If
  \[
  t = p(t) \implies q_1(t) = q_2(t)
  \]
  for all complex numbers $t$, then
  \[
  T \iso p(T) \implies q_1(T) \iso q_2(T)
  \]
  for all objects $T$ of any category in which it makes sense to add and
  multiply objects.
\end{quote}
This is subject to some restrictions on the three polynomials.  For
instance (Theorem~\ref{thm:complex-to-rig}), it suffices to assume that
$p(x) - x$ is primitive, has degree at least two, has non-zero constant
term, and has no repeated complex roots, and that neither $q_1$ nor $q_2$
is constant.  The last condition is what forbids the conclusion $T^6 \iso
1$ in our example.

The story began with one sentence of Lawvere in 1990 (\cite{Lawvere},
p.~11): 
\begin{quote}
  I was surprised to note that an isomorphism $x = 1 + x^2$ [\ldots] always
  induces an isomorphism $x^7 = x$.
\end{quote}
Provoked by this, Blass analysed the situation in detail, producing amongst
other things an explicit bijection between the set $T$ of trees and the set
$T^7$ of 7-tuples of trees; he called the phenomenon `Seven Trees in
One'~\cite{Blass}.  There are many such bijections, none of them
particularly intuitive.  Each corresponds to a way of building an
isomorphism $T \isoarrow T^7$ from a given isomorphism $T \isoarrow 1 +
T^2$ using only multiplication and addition.  One such (not Blass's) runs
as follows: first note that for each $n\geq 1$, we may multiply the given
isomorphism by $T^{n-1}$ (on the left, say) to obtain an isomorphism $T^n
\isoarrow T^{n-1} + T^{n+1}$; use this repeatedly to build a chain of
isomorphisms
\begin{eqnarray*}
T	&\isoarrow	&1 + T^2	\\
	&\isoarrow	&1 + T + T^3	\\
	&\isoarrow	&1 + T + T^2 + T^4	\\
	&\isoarrow	&2T + T^4	\\
	&\vdots		&\vdots		\\
	&\isoarrow	&T^7,
\end{eqnarray*}
with 18 isomorphisms in total.  

We began by trying this method on other polynomials.  For example, we
considered trees in which each vertex has either one or two branches coming
up out of it, the set $T$ of which satisfies $T \iso 1 + T + T^2$.  The
complex solutions are $t = \pm i$, which of course satisfy $t^5 = t$, and
indeed we were able to build an isomorphism $T^5 \iso T$ in a manner
similar to the one for $T^7$ above.  (In fact this example is of special
interest: it leads to a `categorified' or `objectified' version of the
Gaussian integers, as discussed
in~\cite{GINT}.)  More generally, we were able to show that for $n\geq 2$,
\[
T \iso 1 + T + T^n 
\implies
T^{2n+1} \iso T,
\]
and with some effort, found a proof that
\[
T \iso 1 + T + T^2
\implies 
(1 + T)^9 \iso 16(1 + T).
\]
We hoped, of course, that there would turn out to be some general theorem
of which all these isolated results were special cases.  Our hope was
fulfilled, and the subject of this paper is that theorem.  We have
therefore solved the problem posed by Blass in Section~2 of~\cite{Blass}.

Here is the strategy.  Our goal is to turn arguments using complex numbers
into arguments using only addition and multiplication.  Basic commutative
algebra (Section~\ref{sec:rings}) shows that subject to some conditions on
the polynomials involved, if the implication
\[
t = p(t) \implies q_1(t) = q_2(t)
\]
holds for all complex numbers $t$ then it holds for all elements $t$ of all
rings.  We want to conclude that it holds for all rigs (ri\emph{n}gs
without \emph{n}egatives, also known as semirings, Section~\ref{sec:rigs}).
So the challenge is to discover how to turn a proof that uses subtraction
into one that does not.  In precise terms, this is a matter of
cancellability in the underlying additive semigroup of the quotient rig
\[
\nat[x]/(x=p(x)).
\]
We therefore develop (Section~\ref{sec:cancel}) a small amount of
general theory of cancellability in semigroups, and using the assumptions
on $p$, $q_1$ and $q_2$ we establish (Section~\ref{sec:high-polynomials})
the necessary cancellability properties of this particular semigroup.
This is the heart of the paper.

We finish (Section~\ref{sec:main}) by assembling the pieces to give a proof
of the main theorem and looking at some further examples.

\paragraph*{Acknowledgements}  We are grateful to Bill Lawvere and Steve
Schanuel for a remark that gave us the last piece of the jigsaw: see
Section~\ref{sec:cancel}.  We also thank Susan Howson, Peter Johnstone,
Maxim Kontsevich, and Mark Lawson for useful remarks.  

\section{Rings}
\label{sec:rings}

Our rings will always be equipped with multiplicative identities, but need
not be commutative.

\begin{defn}	\label{defn:ring-implication}
Let $p_1, p_2, q_1, q_2 \in \integers[x]$.  We say that
\[
p_1(x) = p_2(x) \implies q_1(x) = q_2(x)
\textup{\demph{ ring-theoretically}}	
\]
if the following equivalent conditions hold:
\begin{enumerate}
\item	\label{item:all-rings}
for all rings $A$ and all $a\in A$, if $p_1(a)=p_2(a)$ then $q_1(a)=q_2(a)$
\item	\label{item:all-comm-rings}
as~\bref{item:all-rings}, but restricted to commutative rings
\item	\label{item:quotient-ring}
$q_1$ and $q_2$ represent the same element of the quotient ring
$\integers[x]/(p_1 - p_2)$ 
\item	\label{item:divisibility}
$(p_1 - p_2)$ divides $(q_1 - q_2)$ in the ring $\integers[x]$.
\end{enumerate}
\end{defn}

(Condition~\bref{item:quotient-ring} implies
condition~\bref{item:all-rings} by the universal property of the quotient,
and the other implications are trivial.)

As suggested by Blass~\cite{Blass}, the first step of `rehabilitating'
nonsense proofs is to rid them of complex numbers and turn them into ring
theory.  So, we start by assuming that
\[
p_1(x) = p_2(x) \implies q_1(x) = q_2(x)
\]
for complex numbers and try to deduce that the same implication holds
ring-theoretically.  This will not work in general: for instance, each of
the implications
\begin{eqnarray}
x = 2 + x + 2x^2	&\implies	&x = 1 + x + x^2	
\label{eq:need-primitive}	\\
x = 1 + 3x + x^2	&\implies	&x = 1 + 2x 
\label{eq:need-separable}	
\end{eqnarray}
holds for complex numbers but fails ring-theoretically.
We therefore seek classes of polynomials for which the deduction is
possible.
\begin{propn}	\label{propn:complex-to-ring}
  Let $p_1, p_2, q_1, q_2 \in \integers[x]$.  Suppose that the polynomial $(p_1
  - p_2) \in \integers[x]$ is primitive and has no repeated complex roots, and
  that each complex root $t$ satisfies $q_1(t) = q_2(t)$.  Then
  \[
  p_1(x) = p_2(x) \implies q_1(x) = q_2(x) \textrm{ ring-theoretically.}
  \]
\end{propn}
(Recall that a polynomial over $\integers$ is \demph{primitive} if the 
greatest common divisor of its coefficients is $1$.
The following proof is little more than Gauss's Lemma: the product of
primitive polynomials is primitive.)
\begin{proof}
By the division algorithm,
\[
q_1 - q_2 = f \cdot (p_1 - p_2) + g
\]
for some $f, g \in \rat[x]$ with $g=0$ or $\deg(g) < \deg(p_1 - p_2)$.
Each of the $\deg(p_1 - p_2)$ complex roots of $p_1 - p_2$ is a root of
$q_1 - q_2$, so of $g$ too; hence $g=0$.  If $f=0$ then certainly
$f\in\integers[x]$, and we are done.  Otherwise we may write $f =
\twid{f}/k$ with $k$ a non-zero integer and $\twid{f} \in \integers[x]$ a
primitive polynomial, and then
\[
k \cdot (q_1 - q_2) = \twid{f} \cdot (p_1 - p_2),
\]
so by Gauss's Lemma $k = \pm 1$ and $f \in \integers[x]$ again.
\done
\end{proof}

Implications~\bref{eq:need-primitive} and~\bref{eq:need-separable} show
that the conditions on primitivity and distinctness of roots cannot be
dropped.

\section{Rigs}
\label{sec:rigs}

A \demph{rig} is a set $A$ equipped with elements $0$ and $1$ and binary
operations $+$ and $\cdot$ such that $(A, 0, +)$ is a commutative monoid,
$(A, 1, \cdot)$ is a monoid, and the distributive laws hold:
\[
\begin{array}{rclcrcl}
0	&=	&a0	&	&0	&=	&0a	\\
ab + ac	&=	&a(b+c)	&	&ba + ca&=	&(b+c)a	
\end{array}
\]
for all $a, b, c \in A$.  

\begin{examples}	\label{egs:rigs}
\begin{enumerate}
\item Any ring is a rig.
\item The initial rig is the set $\nat$ of natural numbers, with its usual
  arithmetic. 
\item The free rig on one generator is the set $\nat[x]$ of polynomials
  over $\nat$. 
\item	\label{eg:rig-lattice}
  Any distributive lattice $A$ is a rig: $+$ is least upper bound, $0$
  is the least element, $\cdot$ is greatest lower bound, and $1$ is the
  greatest element.  

  One might be tempted to try to turn all rigs into rings by adjoining
  negatives formally.  This can certainly be done (and defines a left
  adjoint to the inclusion functor $\fcat{Rings} \arrow \fcat{Rigs}$), but
  destroys a lot of information.  For example, if $A$ is a distributive
  lattice then $a + a = a$ for all $a \in A$, which in the presence of
  negatives implies $a = 0$, so $A$ collapses to the trivial ring.

\item	\label{eg:rig-cardinals}
  The class of all cardinals, with their usual arithmetic, forms a
  large rig.  A set-theoretically respectable version is $\{
  \textrm{cardinals } < \kappa \}$, for any infinite cardinal $\kappa$.

\end{enumerate}
\end{examples}

As with any kind of algebraic structure, it makes sense to talk about
quotients of rigs.  A \demph{congruence} on a rig $A$ is an equivalence
relation $\sim$ on $A$ such that if $\sim$ is regarded as a subset of the
product rig $A \times A$ then it is a sub-rig.  Explicitly, this means
that if $a \sim a'$ and $b \sim b'$ then $a+b \sim a'+b'$ and $ab \sim
a'b'$.  This is precisely the condition needed on the equivalence relation
$\sim$ in order that the set $A/\!\!\sim$ of equivalence classes inherits
the structure of a rig.

\begin{examples}	\label{egs:quotient-rigs}
\begin{enumerate}
\item	\label{eg:rig-of-degrees}
  The relation `has the same degree as' defines a congruence~$\sim$ on
  $\nat[x]$.  We call the quotient $\nat[x]/\!\!\sim$ the \demph{rig of
  degrees}.  It has countably many elements, conveniently written as
  $L^{-\infty}, L^0, L^1, L^2, \ldots$, with operations
  \[
  \begin{array}{rclcrcl}
  L^m + L^n	&=	&L^{\max\{m,n\}}	&	&
  L^m \cdot L^n	&=	&L^{m+n}		\\
  0		&=	&L^{-\infty}		&	&
  1		&=	&L^0
  \end{array}
  \]
  ($m, n \in \{-\infty\} \cup \nat$).  These equations make sense if $L$ is
  thought of as a large number.  

\item	\label{eg:rig-of-codegrees}
  Dually, define the \demph{codegree} of a polynomial $q(x)$ as the least
  $n\in\nat$ for which $q(x)$ has a non-zero coefficient in $x^n$, or as
  $\infty$ if $q=0$.  The quotient of $\nat[x]$ by the congruence `has the
  same codegree as' is the \demph{rig of codegrees}, which has elements
  $\epsln^{\infty}, \epsln^0, \epsln^1, \epsln^2, \ldots$ and operations
  appropriate for $\epsln$ being small.

\end{enumerate}
\end{examples}

Any relation on a rig generates a congruence, which can be defined as the
intersection of all congruences containing the given relation.  We are
particularly interested in the congruence $\sim$ on $\nat[x]$ generated by
declaring equivalent two polynomials $p_1$ and $p_2$.  This can be
described explicitly as follows: $q_1 \sim q_2$ if and only if there is a
finite sequence of polynomials
\[
q_1 = r_0, r_1, \ldots, r_{n-1}, r_n = q_2
\]
(for some $n\in\nat$) such that for each $i \in \{1, \ldots, n\}$, there
exist $f \in \nat[x]$ and $k\in\nat$ satisfying
\begin{equation}	\label{eq:link}
\{ r_{i-1}(x), r_i(x) \}
=
\{ f(x) + x^k p_1(x), f(x) + x^k p_2(x) \}.
\end{equation}
We write the quotient rig $\nat[x]/\!\!\sim$ as $\nat[x]/(p_1 = p_2)$.

The situation for \emph{rings} is much easier: congruences, defined
analogously, correspond to ideals, and if we generate a congruence $\sim$
on the ring $\integers[x]$ by identifying polynomials $p_1$ and $p_2$ then
$q_1 \sim q_2$ if and only if $(p_1 - p_2)$ divides $(q_1 - q_2)$ in
$\integers[x]$.
So the following definition is precisely analogous to
Definition~\ref{defn:ring-implication}.

\begin{defn}	\label{defn:rig-implication}
Let $p_1, p_2, q_1, q_2 \in \nat[x]$.  We say that
\[
p_1(x) = p_2(x) \implies q_1(x) = q_2(x)
\textup{\demph{ rig-theoretically}}	
\]
if the following equivalent conditions hold:
\begin{enumerate}
\item	\label{item:all-rigs}
for all rigs $A$ and all $a\in A$, if $p_1(a)=p_2(a)$ then $q_1(a)=q_2(a)$
\item	\label{item:all-comm-rigs}
as~\bref{item:all-rigs}, but restricted to commutative rigs
\item	\label{item:quotient-rig}
$q_1$ and $q_2$ represent the same element of the quotient rig
$\nat[x]/(p_1 = p_2)$ 
\item	\label{item:chain}
there is a finite sequence $r_0, \ldots, r_n$ of elements of $\nat[x]$
(for some $n\in\nat$) satisfying $r_0 = q_1$, $r_n = q_2$, and the
condition described in~\bref{eq:link}.
\setcounter{bean}{\value{enumi}}
\end{enumerate}
\end{defn}

We want to discuss categories in which objects can be added and multiplied.
Such categories bear the same conceptual relation to rigs as monoidal
categories do to monoids.  So, a \demph{rig category} is a category
$\cat{A}$ equipped with a symmetric monoidal structure $(\oplus, 0)$ and a
monoidal structure $(\otimes, I)$ with the latter distributing over the
former up to coherent isomorphism---in other words, there are specified
isomorphisms
\[
\begin{array}{rclcrcl}
0	&\isoarrow	&T \otimes 0	&	&
0	&\isoarrow	&0 \otimes T	\\
(T \otimes U) \oplus (T \otimes V)	&
\isoarrow	&
T \otimes (U \oplus V)	&	&
(U \otimes T) \oplus (V \otimes T)	&
\isoarrow	&
(U \oplus V) \otimes T	
\end{array}
\]
for each $T, U, V \in \cat{A}$.  The distributivity, associativity and unit
axioms are required to satisfy various axioms; see Laplaza~\cite{Laplaza}
for details.  Any polynomial $p \in \nat[x]$ and object $T$ of a rig
category $\cat{A}$ give rise to a new object $p(T)$ of $\cat{A}$, which the
axioms ensure is well-defined up to canonical isomorphism.

\begin{examples}	
\begin{enumerate}
\item A \demph{distributive category} is a category in which finite
  coproducts and products exist and the latter distribute over the former.
  Any such is naturally a rig category.  Examples are the category of sets,
  the category of topological spaces, any bicartesian closed category, and
  any distributive lattice.
\item The category of sets and partial functions, with disjoint union as
  $\oplus$ and cartesian product as $\otimes$, is a rig category.
\item	\label{eg:real-mens-rig-categories}
  The category of modules over a fixed commutative ring, with the usual
  $\oplus$ and $\otimes$, is a rig category.  The same goes for the
  category of representations of a group and the category of vector bundles
  over a topological space.  
\item A discrete rig category (one in which the only morphisms are the
  identities) is merely a rig.
\end{enumerate}
\end{examples}

The set (or class) of isomorphism classes of objects of a rig category
forms a (possibly large) rig, called its \demph{Burnside rig}.  For
instance, the Burnside rig of the distributive category of sets is the rig
of cardinals~(\ref{egs:rigs}\bref{eg:rig-cardinals}).  The Burnside rigs of
the categories in~\bref{eg:real-mens-rig-categories} are basic to
$K$-theory and representation theory.  (In those subjects the Burnside rig
is no sooner formed than turned into a ring, and, as pointed out
in~\ref{egs:rigs}\bref{eg:rig-lattice}, this process potentially destroys
information: the `Eilenberg swindle' of $K$-theory.
For this reason, among others, the categories
of~\bref{eg:real-mens-rig-categories} are actually replaced by certain
subcategories.)

By considering Burnside rigs and discrete rig categories we see that a
further equivalent condition may be added to
Definition~\ref{defn:rig-implication}:
\begin{enumerate} \itshape
  \setcounter{enumi}{\value{bean}}
  \item for all rig categories $\cat{A}$ and all $T\in \cat{A}$, if $p_1(T)
  \iso p_2(T)$ then $q_1(T) \iso q_2(T)$. 
\end{enumerate}
Moreover, suppose that the rig-theoretic implication holds and that we are
given a specific isomorphism $p_1(T) \isoarrow p_2(T)$, for some $T$ and
$\cat{A}$.  Then there exists a chain $r_0, \ldots, r_n$ of polynomials as
in condition~\ref{defn:rig-implication}\bref{item:chain}, and we can build
from it a specific isomorphism $q_1(T) \isoarrow q_2(T)$.  This is exactly
how the 18-step isomorphism $T \isoarrow T^7$ in the introduction was
built. 

If a proof can be done without using subtraction then it can certainly be
done with subtraction available; in other words, if an implication holds
rig-theoretically then it certainly holds ring-theoretically.  This paper
is about going the other way, and the next result shows that it is a
question of cancellability.
\begin{propn}	\label{propn:ring-to-rig-plus}
  Let $p_1, p_2, q_1, q_2 \in \nat[x]$ and suppose that
  \[
  p_1(x) = p_2(x) \implies q_1(x) = q_2(x)
  \textrm{ ring-theoretically}.
  \]
  Then there exists $s \in \nat[x]$ such that 
  \[
  p_1(x) = p_2(x) \implies q_1(x) + s(x) = q_2(x) + s(x)
  \textrm{ rig-theoretically}.
  \]
\end{propn}
\begin{proof}
We are given that there exists $r \in \integers[x]$ satisfying
\[
q_1 - q_2 = r \cdot (p_1 - p_2)
\]
in $\integers[x]$.  We may write $r = r_1 - r_2$ for some $r_1, r_2 \in
\nat[x]$,  and then
\[
q_1 + r_1 p_2 + r_2 p_1 = q_2 + r_1 p_1 + r_2 p_2
\]
in $\nat[x]$.  Put $s = r_1 p_1 + r_2 p_2$: then $q_1 + s$ and $q_2
+ s$ represent the same element of the quotient rig $\nat[x]/(p_1 = p_2)$,
as required.
\done
\end{proof}

\section{Cancellation in commutative semigroups}
\label{sec:cancel}

A \demph{commutative semigroup} $(A, \csg)$ is a set $A$ equipped with a
commutative associative binary operation $\csg$.  In general, $a_1 \csg b =
a_2 \csg b$ does not imply $a_1 = a_2$, but in this section we give a
condition under which it does.  

Later we will apply this to the underlying additive semigroup of a rig, but
it seems to be easier to understand the following results if the semigroup
operation $\csg$ is thought of as multiplication.  Informally, take a
commutative semigroup and call an element `high' if every element divides
it.  Then the set of high elements is closed under multiplication, and in
it every element divides every other element.  This says that the set of
high elements is, if not empty, a group.  So given an equation $a_1 \csg b
= a_2 \csg b$ in which each $a_i$ is high, we may post-multiply each side
by $a_1$ then divide through by the high element $b \csg a_1$ to conclude
that $a_1 = a_2$.

Formally, given a commutative semigroup $(A, \csg)$, define a relation
$\leq_A$ on $A$ by
\[
b \leq_A a 
\iff
\textrm{there exists } c\in A \textrm{ satisfying } b \csg c = a.
\]
The notation is potentially misleading: $\leq_A$ is transitive but not
necessarily reflexive (consider strictly positive numbers under addition)
or antisymmetric (consider an abelian group).  However, $\leq_A$ has the
expected meaning when $(A, \csg) = (\nat, +)$ or when $\csg$ is the least
upper bound operation on a (semi-)lattice
(\emph{cf}.~\ref{egs:rigs}\bref{eg:rig-lattice}).

An element $a$ of $A$ is called \demph{high} if $b \leq_A a$ for all $b \in
A$, and the set of high elements of $A$ is written $\High(A)$.  This set
may be empty (as for $(\nat,+)$), or all of $A$ (as for abelian groups), or
somewhere in between (interesting examples of which occur later).  We call
$A$ a \demph{clique} if $b \leq_A a$ for all $a, b \in A$, or,
equivalently, if $\High(A) = A$.

\begin{lemma}
  Let $(A, \csg)$ be a commutative semigroup.  Then $\High(A)$ is a
  sub-semigroup of $A$, and $(\High(A), \csg)$ is a clique.  
\end{lemma}
\begin{proof}
Let $a, b \in \High(A)$.  We have to show that $a\csg b \in \High(A)$ and
that there exists $c \in \High(A)$ satisfying $a = b \csg c$.  In fact, we
can do both without knowing that $b$ is high.  First, we have $a \csg b
\geq_A a \in \High(A)$, hence $a \csg b \in \High(A)$.  Second, we have $a
\geq_A b \csg a$, so there exists $d \in A$ satisfying $a = b \csg a \csg
d$; take $c = a \csg d \in \High(A)$.
\hspace*{80mm}\done
\end{proof}

We observed earlier that any abelian group is a clique.  The next result is
very nearly the converse.
\begin{lemma}	\label{lemma:cliques}
  A commutative semigroup is a clique if and only if it is empty or an
  abelian group.
\end{lemma}
\begin{proof}
Let $(B, \csg)$ be a nonempty clique, and pick an element $d$.  Since $d
\leq_B d$, there exists $\hz \in B$ satisfying $d \csg \hz = d$.
For any $b\in B$ we have $b \geq_B d$, so there exists $c\in B$ satisfying
$b = d \csg c$, and then using commutativity, $b \csg \hz = b$.
Hence $\hz$ is a unit for $\csg$.  Also, $B$ has inverses: if $b
\in B$ then $b \leq_B \hz$. 
\done
\end{proof}

\begin{cor}	\label{cor:abelian-group}
  If $A$ is a commutative semigroup then $\High(A)$ is either empty or, when
  equipped with the inherited binary operation, an abelian group.  \done
\end{cor}

This corollary may be in the semigroup literature, 
but we have been unable to find it.  

Beware that even if the semigroup $A$ has a unit, this is very likely
\emph{not} the unit of $\High(A)$.  (Indeed, the units can only be the same
if $A$ is an abelian group to start with.)  For a simple example, take $A$
to be a lattice and $\csg$ to be least upper bound: then the unit of $A$ is
the least element, but $\High(A)$ is the singleton consisting of the
greatest element.  

\begin{cor}	\label{cor:cancellation}
  Let $a_1$ and $a_2$ be high elements of a commutative semigroup $(A,
  \csg)$.  If there exists $b\in A$ such that $a_1 \csg b = a_2 \csg b$
  then $a_1 = a_2$.  \done
\end{cor}

As explained, we will apply these results when $A$ is the underlying
additive semigroup of a rig.  In that case, although we will not need to
know it, we have:
\begin{cor}	\label{cor:high-ring}	
  If $A$ is a rig then $\High(A,+)$ is either empty or, when equipped with
  the inherited binary operations, a ring.
\end{cor}
\begin{proof}
In the notation of~\ref{lemma:cliques}, the multiplicative unit is $1 +
\hz$.  The only non-trivial check is that $a \cdot \hz = \hz$ for all high
$a$.  \done
\end{proof}

We arrived at the results of this section after a conversation between one
of us (M.F.) and Bill Lawvere, in which Lawvere mentioned a result of Steve
Schanuel's that `the infinite-dimensional elements form a ring'.  But we do
not know any details of this work beyond what is in~\cite{Schanuel}.

\section{High polynomials}
\label{sec:high-polynomials}

Our final task is to find the high elements of the quotient rig
$\nat[x]/(p_1 = p_2)$ (or rather, its underlying additive semigroup).  We
do this under several assumptions, the most sweeping of which is that
$p_1(x) = x$.  

Fix a polynomial $p(x) \in \nat[x]$. 
To lighten the notation, we say that a polynomial $f$ is high to mean that
its image $[f]$ in the quotient $\nat[x]/(x = p(x))$ is high, and write $g
\leq f$ to mean $[g] \leq_{\nat[x]/(x = p(x))} [f]$.  Observe that $\leq$
is compatible with multiplication: if $g_1 \leq f_1$ and $g_2 \leq f_2$
then $g_1 g_2 \leq f_1 f_2$.

\begin{lemma}	\label{lemma:ladder}
  If $p$ has non-zero constant term then $1 \leq x \leq x^2 \leq \cdots$.
\end{lemma}
\begin{proof}
  We have $p(x) \geq 1$ by hypothesis, so $x\geq 1$, and multiplying
  through by $x^n$ gives $x^{n+1} \geq x^n$ for all $n\in\nat$.  \done
\end{proof}

\begin{lemma}   \label{lemma:multiples+powers}
If $p$ has non-zero constant term and degree at least two then 
\begin{enumerate}
\item \label{item:multiples}
  $x \geq nx$ for all $n\in\nat$, and
\item \label{item:powers}
  $x \geq x^n$ for all $n\in\nat$. 
\end{enumerate}
\end{lemma}
\begin{proof}
The hypotheses imply that $x \geq 1 + x^d$ for some $d\geq 2$.  By
induction, it is enough to prove each of~\bref{item:multiples}
and~\bref{item:powers} when $n = 2$; and using Lemma~\ref{lemma:ladder},
\[
x \geq x^d \geq x^{d-1} (1 + x^d) = x^{d-1} + x^{2d-1} \geq 2x
\]
and
\[
x \geq x^d \geq x^2.
\]
\ \done
\end{proof}

\begin{propn}   \label{propn:non-constant}
If $p$ has non-zero constant term and degree at least two then
every non-constant polynomial is high.
\end{propn}
\begin{proof}
By Lemma~\ref{lemma:multiples+powers},
\[
x \ \geq \ kx 
  \ \geq \ x^{n_1} + x^{n_2} + \cdots + x^{n_k}
\]
for all $k, n_1, \ldots, n_k \in \nat$; in other words, $x$ is high.  The
result then follows from Lemma~\ref{lemma:ladder}.
\done
\end{proof}

We will not need to know it, but under the hypotheses of the Proposition,
the high polynomials are precisely the non-constants.  This can
be proved by considering the 3-element quotient rig of $\nat[x]$ in which
the equivalence classes are $\{0\}$, the set of non-zero constants, and the
set of non-constants.  So by Corollary~\ref{cor:high-ring}, the
non-constant polynomials form a ring.  The quotient rig $\nat[x]/(x=p(x))$
is the disjoint union of this ring with the set of natural numbers.

\section{The Main Theorem}
\label{sec:main}

We have now done all the work and can read off the main theorem, of
which we give two slightly different versions.

\begin{thm}	\label{thm:ring-to-rig}
  Let $p, q_1, q_2 \in \nat[x]$ be polynomials such that $p$ has non-zero
  constant term and degree at least two and $q_1$ and $q_2$ have degree at
  least one.  If
  \[
  x = p(x) \implies q_1(x) = q_2(x) \textrm{ ring-theoretically}
  \]
  then the same is true rig-theoretically.
\end{thm}
\begin{proof}
  Assemble Proposition~\ref{propn:ring-to-rig-plus},
  Corollary~\ref{cor:cancellation}, and
  Proposition~\ref{propn:non-constant}.  \done
\end{proof}

\begin{thm}	\label{thm:complex-to-rig}
  Let $p, q_1, q_2$ be polynomials as in the first sentence of
  Theorem~\ref{thm:ring-to-rig}.  Suppose that the polynomial $p(x) -
  x \in \integers[x]$ is primitive and has no repeated complex roots, and
  that each complex root $t$ satisfies $q_1(t) = q_2(t)$.  Then
  \[
  x = p(x) \implies q_1(x) = q_2(x) \textrm{ rig-theoretically.}
  \]
\end{thm}
\begin{proof}
  Assemble Proposition~\ref{propn:complex-to-ring}
  and Theorem~\ref{thm:ring-to-rig}.  \done
\end{proof}

The first version is more general, and is the one to use when the complex
solutions of $x = p(x)$ are hard to find.
To apply it we verify by the division algorithm that the ring-theoretic
implication holds.  The second version is more useful when the solutions of
$x = p(x)$ are known or easily calculated, and is applied by simply
checking that $q_1(t) = q_2(t)$ for each solution $t$.

The proofs are constructive.  So, if we are given polynomials $p, q_1, q_2$
satisfying the hypotheses of either~\ref{thm:ring-to-rig}
or~\ref{thm:complex-to-rig} then unwinding the proof gives an explicit
sequence of polynomials
demonstrating the rig-theoretic implication (as in
Definition~\ref{defn:rig-implication}\bref{item:chain}).  Combining this
with the observations of Section~\ref{sec:rigs}, if we are also given an
isomorphism $T \isoarrow p(T)$ for some object $T$ of a rig category then
we obtain an explicit isomorphism $q_1(T) \isoarrow q_2(T)$.

\begin{examples}
\begin{enumerate}

\item  Returning to our original example, $1 - x + x^2$ has distinct
  complex roots $e^{\pm i\pi/3}$, both of which satisfy $x^6 = 1$, so
  Theorem~\ref{thm:complex-to-rig} tells us that
  \[
  T \iso 1 + T^2 \implies T^7 \iso T
  \]
  for any object $T$ of a rig category.  (Henceforth we write $+$ and
  $\times$ for the monoidal structures of a rig category.)  In particular,
  $T^7 \iso T$ when $T$ is the set of binary trees.  

\item If $n\geq 2$ then the solutions of the equation $x = 1 + x + x^n$ are
  the complex $n$th roots of $-1$, each of which satisfies $x^{2n} = 1$.
  Hence
  \[
  T \iso 1 + T + T^n \implies T^{2n+1} \iso T 
  \textrm{ and }
  T + T^{2n} = 1 + T
  \]
  for any object $T$ of a rig category.

\item The case $n=2$ of the previous example is particularly easy to work
  with since the solutions are $\pm i$.  For instance, if $T$ is an object
  of a rig category satisfying $T \iso 1 + T + T^2$ then there are also
  isomorphisms
  \[
  T^4 \iso 2 + T^2,
  \diagspace  
  T + T^3 \iso 1 + T^2,
  \diagspace
  (1 + T)^9 \iso 16(1 + T)
  \]
  (observing for the last one that $1 \pm i = \sqrt{2} e^{\pm i\pi/4}$).
  Our paper~\cite{GINT} explores this example in depth.

\item 
Let $m$ and $n$ be coprime positive integers, one of which is
even.  Then the complex roots of $(1 + x^m)(1 + x^n)$ 
are distinct and each satisfy $x^{2mn} = 1$, so
\[
T \iso 1 + T + T^m + T^n + T^{m+n}
\]
implies
\[
  T^{2mn + 1} \iso T 
  \textrm{ and } T + T^{2mn} \iso 1 + T
\]
for any object $T$ of a rig category.

\item 
The following randomly-chosen example illustrates the power and generality of
Theorem~\ref{thm:ring-to-rig}.  Let
\begin{eqnarray*}
p(x)	&=	&3 + 2x^3 + 4x^5,	\\
q_1(x)	&=	&6x + 10x^2 + x^3 + 3x^4 +2x^5 + 7x^6 + 12x^7,	\\
q_2(x)	&=	&3 + 2x + 2x^2 + 9x^3 + 5x^6 + 4x^8.
\end{eqnarray*}
A routine application of the division algorithm shows that $p(x) - x$
divides $q_1(x) - q_2(x)$ in $\integers[x]$, so by
Theorem~\ref{thm:ring-to-rig},
\[
x = p(x) \implies q_1(x) = q_2(x)
\]
rig-theoretically.  Certainly this implication would be tiresome to prove
by hand.  Indeed, without the results of this paper it would not be at all
clear that there was any systematic way of finding such a proof.

\end{enumerate}
\end{examples}

Observe finally that Theorems~\ref{thm:ring-to-rig} and
Theorem~\ref{thm:complex-to-rig} are sharp: none of the hypotheses can be
dropped. 
For the condition that $p$ has non-zero constant
term, consider the implication
\[
x = x + x^2 \implies x^2 = x^3.
\]
This holds ring-theoretically, but fails when $x$ is the element $\epsln^1$
of the rig of
codegrees~(\ref{egs:quotient-rigs}\bref{eg:rig-of-codegrees}).  For the
condition that $p$ has degree at least two, consider
\[
x = 1 + x \implies x = x^2,
\]
which holds ring-theoretically but fails when $x$ is the element $L^1$ of
the rig of degrees~(\ref{egs:quotient-rigs}\bref{eg:rig-of-degrees}).  For
the condition that $q_1$ and $q_2$ are non-constant, consider the original
example of
\[
x = 1 + x^2 \implies x^6 = 1,
\]
which holds ring-theoretically but fails when $x$ is the element $\aleph_0$
of the rig of countable cardinals.  And we saw in Section~\ref{sec:rings}
that the extra hypotheses in Theorem~\ref{thm:complex-to-rig} (primitivity
and distinctness of roots) cannot be dropped, otherwise the implication
might not even hold ring-theoretically.

\end{document}